# A Note on Analyzing the Stability of Oscillator Ising Machines

Mohammad Khairul Bashar, Zongli Lin, and Nikhil Shukla*

*Department of Electrical and Computer Engineering, University of Virginia, Charlottesville, VA- 22904 USA*

*e-mail: ns6pf@virginia.edu

The rich non-linear dynamics of the coupled oscillators (under second harmonic injection) can be leveraged to solve computationally hard problems in combinatorial optimization such as finding the ground state of the Ising Hamiltonian. While prior work on the stability of the so-called Oscillator Ising Machines (OIMs) has used the linearization method, in this letter, we present a complementary method to analyze stability using the second order derivative test of the energy / cost function. We establish the equivalence between the two methods, thus augmenting the tool kit for the design and implementation of OIMs.

*Introduction:* Ising Machines (OIMs) offer a markedly different paradigm for solving combinatorial optimization – a problem class that finds increasing importance in areas from machine learning (e.g., training of Boltzmann machines [1]) to applications in logistics such as route planning [2]. Despite the ample use cases of such problems, solving them efficiently using digital computers and algorithms remains a challenge [3]. Consequently, this has motivated the quest for alternate computing paradigms such as Ising Machines to accelerate the solution of such problems [4].

Networks of coupled oscillators offer a promising hardware platform to realize Ising Machines, commonly referred to as, Oscillator Ising Machines (OIMs) [5]. Such a platform offers the advantage of being compact, scalable, low power, and compatible with commercial CMOS foundry process technology [6]. The underlying design principle of such a platform exploits the natural equivalence between the global minima of the energy / cost function ($E$) describing a network of coupled oscillators under second harmonic injection and the global minima of the Ising Hamiltonian ($H = -\sum_{i,j,i<j}^{N} W_{ij}s_i s_j$, where Ising spin $s_i \in \{-1,+1\}$, and $W_{ij} = -E_{ij}$, with $E_{ij}$ being the edge weight between node $i$ and $j$) [7]. The dynamics and the corresponding cost function of such an OIM can be respectively described by:

$$\frac{d\theta_i}{dt} = -K \sum_{i=1, i<j}^{N} W_{ij}\sin(\theta_i - \theta_j) - K_s \sin(2\theta_i(t)) = f_i(\theta) \quad (1)$$

$$= -\frac{1}{2}(\nabla E)_i$$

$$E(\theta(t)) = -K \sum_{i,j=1, j\neq i}^{N} W_{ij}\cos(\theta_i - \theta_j) \quad (2)$$
$$- K_s \sum_{i=1}^{N} \cos(2\theta_i(t))$$

Where, $\theta = (\theta_1, \theta_2, \theta_3, ..., \theta_N)$ represents the oscillator phases, and $K$ and $K_s$ are the coupling and second harmonic injection strength, respectively. From equation (1) and (2), it can be deduced that,

$$\frac{dE}{dt} = \sum_{i=1}^{N} \frac{\partial E}{\partial \theta_i}\frac{d\theta_i}{dt} = -2\sum_{i=1}^{N} \left(\frac{d\theta_i}{dt}\right)^2 \leq 0 \quad (3)$$

The dynamics described in equation (1) reveal that the system of equations has multiple fixed points, where $\frac{d\theta}{dt} = 0$ (also, $-(\nabla E) = 0$). In fact, every spin configuration $(\theta_1, \theta_2, \theta_3, ..., \theta_N)$, $\theta_i \in \{0, \pi\}$ is a fixed point, resulting in $2^N$ such fixed points. We note that there might be other fixed points defined by $(\theta_1, \theta_2, \theta_3, ..., \theta_N)$ where $\theta_i \notin \{0, \pi\}$. The fixed points $(\theta_1, \theta_2, \theta_3, ..., \theta_N)$ with the lowest energy represent the optimal solutions to the Ising Hamiltonian while the rest represent sub-optimal solutions. In our prior work [8], we analyzed the stability of the fixed points, corresponding to various spin configurations, of the OIM. Specifically, we showed that tuning the ratio of the coupling strength among the oscillators ($K$) and the strength of the second harmonic injection ($K_s$) has a dramatic impact on the stability of the globally optimal and sub-optimal fixed points. This analysis was performed by using the linearization method, i.e., by computing the Lyapunov exponents, which in this case, are defined by the eigenvalues of the Jacobian matrix ($J$). A fixed point (globally optimal or not) is attractive if all eigenvalues at that point are negative; and unstable if at least one eigenvalue is positive.

*Results:* While the Jacobian analysis essentially entails working with the first order derivatives, the purpose of this work is to present an alternate approach, based on the second order derivates test of $E$ (using the Hessian Matrix). We show that for a system whose dynamics are of the form $-\alpha(\nabla E)_i = \frac{d\theta_i}{dt}$ ($\alpha = \frac{1}{2}$ for OIM), the stability of the fixed points can be analyzed using the eigenvalues of the Hessian Matrix ($H_E$). A fixed point is attractive if all eigenvalues of $H_E$ are positive. For the unstable case, if some or all of the eigenvalues are negative, then the fixed point is a saddle point, or a local maximum, respectively. To prove this, we establish the equivalence between the two methods.

The Hessian matrix for $E(\theta)$ can be computed as,

$$H_E = \begin{bmatrix} \frac{\partial^2 E}{\partial\theta_1\partial\theta_1} & \frac{\partial^2 E}{\partial\theta_1\partial\theta_2} & \cdots & \frac{\partial^2 E}{\partial\theta_1\partial\theta_N} \\ \frac{\partial^2 E}{\partial\theta_2\partial\theta_1} & \frac{\partial^2 E}{\partial\theta_2\partial\theta_2} & \cdots & \frac{\partial^2 E}{\partial\theta_2\partial\theta_N} \\ \vdots & \vdots & \ddots & \vdots \\ \frac{\partial^2 E}{\partial\theta_N\partial\theta_1} & \frac{\partial^2 E}{\partial\theta_N\partial\theta_2} & \cdots & \frac{\partial^2 E}{\partial\theta_N\partial\theta_N} \end{bmatrix}$$

$$= \begin{bmatrix} \frac{\partial^2 E}{\partial\theta_1\partial\theta_1} & \frac{\partial^2 E}{\partial\theta_2\partial\theta_1} & \cdots & \frac{\partial^2 E}{\partial\theta_N\partial\theta_1} \\ \frac{\partial^2 E}{\partial\theta_1\partial\theta_2} & \frac{\partial^2 E}{\partial\theta_2\partial\theta_2} & \cdots & \frac{\partial^2 E}{\partial\theta_N\partial\theta_2} \\ \vdots & \vdots & \ddots & \vdots \\ \frac{\partial^2 E}{\partial\theta_1\partial\theta_N} & \frac{\partial^2 E}{\partial\theta_2\partial\theta_N} & \cdots & \frac{\partial^2 E}{\partial\theta_N\partial\theta_N} \end{bmatrix} \quad (4)$$

Here, the symmetricity of the second derivatives ($\frac{\partial^2 E}{\partial\theta_i\partial\theta_j} = \frac{\partial^2 E}{\partial\theta_j\partial\theta_i}$) is used. Using Eq. (1), $\frac{d\theta_i}{dt} = -\frac{1}{2}(\nabla E)_i = -\frac{1}{2}\frac{\partial E}{\partial\theta_i} = f_i$, equation (4) can now be expressed as:

$$H_E = \begin{bmatrix} \frac{\partial^2 E}{\partial\theta_1\partial\theta_1} & \frac{\partial^2 E}{\partial\theta_2\partial\theta_1} & \cdots & \frac{\partial^2 E}{\partial\theta_N\partial\theta_1} \\ \frac{\partial^2 E}{\partial\theta_1\partial\theta_2} & \frac{\partial^2 E}{\partial\theta_2\partial\theta_2} & \cdots & \frac{\partial^2 E}{\partial\theta_N\partial\theta_2} \\ \vdots & \vdots & \ddots & \vdots \\ \frac{\partial^2 E}{\partial\theta_1\partial\theta_N} & \frac{\partial^2 E}{\partial\theta_2\partial\theta_N} & \cdots & \frac{\partial^2 E}{\partial\theta_N\partial\theta_N} \end{bmatrix}$$

$$= \begin{bmatrix} \frac{\partial(-2f_1)}{\partial\theta_1} & \frac{\partial(-2f_1)}{\partial\theta_2} & \cdots & \frac{\partial(-2f_1)}{\partial\theta_N} \\ \frac{\partial(-2f_2)}{\partial\theta_1} & \frac{\partial(-2f_2)}{\partial\theta_2} & \cdots & \frac{\partial(-2f_2)}{\partial\theta_N} \\ \vdots & \vdots & \ddots & \vdots \\ \frac{\partial(-2f_N)}{\partial\theta_1} & \frac{\partial(-2f_N)}{\partial\theta_2} & \cdots & \frac{\partial(-2f_N)}{\partial\theta_N} \end{bmatrix}$$

$$= -2\begin{bmatrix} \frac{\partial f_1}{\partial\theta_1} & \frac{\partial f_1}{\partial\theta_2} & \cdots & \frac{\partial f_1}{\partial\theta_N} \\ \frac{\partial f_2}{\partial\theta_1} & \frac{\partial f_2}{\partial\theta_2} & \cdots & \frac{\partial f_2}{\partial\theta_N} \\ \vdots & \vdots & \ddots & \vdots \\ \frac{\partial f_N}{\partial\theta_1} & \frac{\partial f_N}{\partial\theta_2} & \cdots & \frac{\partial f_N}{\partial\theta_N} \end{bmatrix} = -2J \quad (5)$$

Equation (5) reveals that in OIMs, the Jacobian matrix is half of the negative of the Hessian matrix ($J = -\frac{1}{2}H_E$). Consequently, the magnitude of the eigenvalues of the Jacobian matrix ($\lambda_J$) will be half of the magnitude of the eigenvalues of the Hessian matrix ($\lambda_{H_E}$) but with the opposite sign, i.e., $\lambda_J = -\frac{1}{2}\lambda_{H_E}$. This implies that the condition for a fixed point to be attractive, when analyzing the Hessian matrix, is that all eigenvalues be positive. When some eigenvalues are negative, it entails that the fixed point is a saddle point; when all eigenvalues are negative, it can be inferred that the fixed point is a maximum. For the general class of gradient descent systems defined by $-\alpha(\nabla E)_i = \frac{d\theta_i}{dt}$ ($\alpha > 0$), $\lambda_J = -\alpha\lambda_{H_E}$.

*Conclusion:* In summary, we have presented a complementary method that utilizes the eigenvalues of the Hessian matrix to compute the stability of the fixed points in the OIM and is generally applicable to the broader class of gradient descent systems whose dynamics can be defined by $-\alpha(\nabla E)_i = \frac{d\theta_i}{dt}$. Our approach shows that the stability of the fixed points in the OIM can be analyzed by identifying whether they represent a minimum of the energy function.

*Acknowledgments:* This work was supported by NSF grant 2328961.

*Competing Interests:* The authors declare no competing interests.

*Data Availability:* The data that support the findings of this study are available from the corresponding author upon reasonable request.

*Author Contributions:* Mohammad Khairul Bashar: Conceptualization (lead), Writing – Original Draft (supporting). Zongli Lin: Validation (lead), Writing – Review & Editing (lead). Nikhil Shukla: Conceptualization (Supporting), Supervision (lead), Project administration (lead), Writing – Original Draft (lead).

**References**

1. Niazi, S. Aadit, N.A. Mohseni, M. Chowdhury, S. Qin, Y. & Camsari, K.Y.: Training Deep Boltzmann Networks with Sparse Ising Machines. *arXiv preprint* arXiv:2303.10728 (2023).
2. Bao, S. Tawada, M. Tanaka, S. & Togawa, N.: Multi-day travel planning using Ising machines for real-world applications. In *2021 IEEE International Intelligent Transportation Systems Conference (ITSC)*, 3704-3709 (IEEE, 2021).
3. Calude, C.S.: Unconventional computing: A brief subjective history. *Springer International Publishing* (2017).
4. Mohseni, N. et al.: Ising machines as hardware solvers of combinatorial optimization problems. *Nature Reviews Physics* **4**, 363-379 (2022).
5. Bashar, M.K. et al.: Experimental demonstration of a reconfigurable coupled oscillator platform to solve the Max-cut problem. *IEEE JXCDC* **6**, pp. 116-121 (2020).
6. Mallick, A. et al.: Overcoming the accuracy vs. performance trade-off in oscillator ising machines. In *2021 IEDM*, 40-2 (IEEE, 2021).
7. Wang, T. et al.: Solving combinatorial optimisation problems using oscillator based Ising machines. *Natural Computing* **20**, 287-306 (2021).
8. Bashar, M.K. Lin, Z. & Shukla, N.: Stability of Oscillator Ising Machines: Not All Solutions Are Created Equal. *Journal of Applied Physics* **134**, 144901 (2023).